\documentclass[12pt]{article}

\usepackage[body={16cm,23cm}]{geometry}
\usepackage{epsfig,amsfonts,amssymb}
\usepackage{enumerate,amsmath}
\usepackage{cite}
\usepackage[mathscr]{eucal}
\newcommand{\bZ}{{\mathbb Z}}
\def\p{\partial}
\def\a{\alpha}
\def\b{\beta}
\def\g{\gamma}
\def\d{\delta}

\def\beq{\begin{equation}}
\def\eeq{\end{equation}}
\def\be{\begin{displaymath}}
\def\ee{\end{displaymath}}
\def\bea{\begin{eqnarray}}
\def\eea{\end{eqnarray}}

\def\bmat{\left(\begin{array}}

\def\ds{\displaystyle}

\def\Wp{ \raise.4ex\hbox{\textrm{\Large $\wp$}}}

\newtheorem{theorem}{Theorem}[section]

        \newcommand{\qed}{\nobreak \ifvmode \relax \else
      \ifdim\lastskip<1.5em \hskip-\lastskip
      \hskip1.5em plus0em minus0.5em \fi \nobreak
      \vrule height0.75em width0.5em depth0.25em\fi}%%%%%%%%%%%%%%%%%%%%%%%%%%%%%%%%%%%%%%%%%%
\def\eqref#1{(\ref{#1})}
\def\?{(?)\marginpar{|?}}

\begin{document}

\title{Picard solution of Painlev\'e VI and related tau-functions}

\author{
                Vladimir V. Mangazeev\footnote{email:
                {\tt Vladimir.Mangazeev@anu.edu.au}}\\ \\
 Department of Theoretical Physics,\\
         Research School of Physics and Engineering,\\
    Australian National University, \\
    Canberra, ACT 0200, Australia.}
%\date{24 November, 2009}
\maketitle

\begin{abstract}
In this paper we obtain explicit expressions for tau-functions
related to Picard type solutions of the Painlev\'e VI equation
in terms of theta functions and their derivatives.
\end{abstract}

%%%%%%%%%%%%%%%%%%%%%%%%%%%%%%%%%%%%%%%%%%%%%%%%%%

\section{Introduction}

In this paper we study a special case of the Painlev\'e VI equation
\cite{Pain02,Gamb10}
\bea
&{\ds
q''(t)=\frac{1}{2}\biggl(\frac{1}{q(t)}
+\frac{1}{q(t)-1}+\frac{1}{q(t)-t}\biggr)q'(t)^2-
\biggl(\frac{1}{t}+\frac{1}{(t-1)}+\frac{1}{q(t)-t}\biggr)q'(t)+}&
\nonumber\\
&{\ds+
\frac{q(t)(q(t)-1)(q(t)-t)}{t^2(t-1)^2}\biggl[\a+\b\frac{t}{q(t)^2}+
\g\frac{t-1}{(q(t)-1)^2}+\d\frac{t(t-1)}{(q(t)-t)^2}\biggr]
}&\ ,\label{pain1}
\eea
when
\beq
\a=0,\quad \b=0,\quad \g=0,\quad \delta=\frac{1}{2}.\label{pain2}
\eeq

This case was originally considered by Picard \cite{Pic89}.
Due to a special choice of parameters (\ref{pain2}) a general solution
of (\ref{pain1}) is known
\beq
q_0(t)=\textrm{\Large$\wp$}(c_1\omega_1+c_2\omega_2;\omega_1,\omega_2)
+\frac{t+1}{3},\label{pain3}
\eeq
where {\Large$\wp$}$(u;\omega_1,\omega_2)$ is the Weierstrass elliptic
 function
with half-periods $\omega_{1,2}$, $c_{1,2}$ are complex constants
 and $\omega_{1,2}(t)$ are two linearly independent solutions
of the hypergeometric equations
\beq
t(1-t)\omega''(t)+(1-2t)\omega'(t)-\frac{1}{4}\omega(t)=0.\label{pain4}
\eeq

The properties of the Picard solutions have been studied recently
by M. Mazzocco \cite{Maz01}. In particular, she investigated
its monodromy properties and
algebraic solutions
which correspond to $c_1$ and $c_2$ being rational numbers.

Algebraic solutions of Painlev\'e VI play an important role
in many applications in theoretical physics (see, for example,
\cite{Dub96,BM}). In such cases a calculation of related
tau-functions
 can be
simpler due to a presence of an algebraic relation between
the solution $q(t)$
and the variable $t$.

The goal of this paper is different.
We aim to present
explicit expressions for tau-functions
related to the Picard solutions and its images
under birational canonical transformations
\cite{Ok87}
 for {\it generic} values of complex parameters $c_1$ and $c_2$.
To our knowledge this has not been done before.

More explicitly, we consider a sequence of tau-functions \cite{Ok87}
obtained by a parallel shift $l_3$ from the Picard solution.
We calculate
the first two tau-functions and other members of the sequence
 can be obtained using the standard Toda-type
second order relations.

\section{Properties of the Painlev\'e VI equation}

In this section we briefly review the main properties of the equation
(\ref{pain1}) which we denote as
${\bf P}_{ VI}(\a,\b,\g,\d)$.

Following \cite{Ok87} one can introduce two different
parameterizations of parameters in the Painlev\'e VI equation:
$(\kappa_0,\kappa_1,
\kappa_\infty,\theta)$
\beq
\a=\frac{1}{2}\kappa_\infty^2,\quad
\b=-\frac{1}{2}\kappa_0^2,\quad
\g=\frac{1}{2}\kappa_1^2,\quad
\d=\frac{1}{2}(1-\theta^2),\label{pain5}
\eeq
and $(b_1,b_2,b_3,b_4)$
\beq
\kappa_0=b_1+b_2,\quad \kappa_1=b_1-b_2,
\quad \kappa_\infty=b_3-b_4,\quad \theta=b_3+b_4+1.\label{pain6}
\eeq

This equation is equivalent to the Hamiltonian system $H_{VI}(t;q,p)$
described by the equations
\beq
\frac{dq}{dt}=\frac{\partial H}{\partial p},
\quad \frac{dp}{dt}=-\frac{\partial H}{\partial q}\ ,\label{pain7}
\eeq
with
the Hamiltonian function
\bea
&{\ds H_{VI}(t;q,p)=\frac{1}{t(t-1)}\Bigl[q(q-1)(q-t)p^2-}&\nonumber\\
&{\ds -\{\kappa_0\,(q-1)(q-t)+
\kappa_1\,q(q-t)+(\theta-1)q(q-1)\}p+\kappa\>(q-t)],}&\label{pain8}
\eea
where $q\equiv q(t)$, $p\equiv p(t)$ and
\beq
\kappa=\frac{1}{4}(\kappa_0+\kappa_1+\theta-1)^2
-\frac{1}{4}\kappa_\infty^2.\label{pain9}
\eeq
One can introduce an auxiliary Hamiltonian $h(t)$,
\beq
h(t)=t(t-1)H(t)+e_2(b_1,b_3,b_4)
\,t-\frac{1}{2}e_2(b_1,b_2,b_3,b_4),\label{pain10}
\eeq
where $e_i(x_1,\ldots,x_n)$ is the $i$-th elementary symmetric
function in $n$ variables
and a set of $x_i$'s can be a subset of $b_i$'s as in (\ref{pain10}).

Okamoto \cite{Ok87} showed that for each pair $\{q(t),p(t)\}$
satisfying (\ref{pain7}), the function $h(t)$ solves
the  ${\bf E}_{{VI}}$ equation which is
\beq
h'(t)\Bigl[t(1-t)h''(t)\Bigr]^2
+\Bigl[h'(t)[2h(t)-(2t-1)h'(t)]+b_1b_2b_3b_4\Bigr]^2=
\prod_{k=1}^4\Bigl(h'(t)+b_k^2\Bigr) \label{pain11}
\eeq
and $q(t)$ solves ${\bf P}_{ VI}(\a,\b,\g,\d)$.

Conversely, for each solution $h(t)$ of (\ref{pain11}),
such that $\frac{d^2}{dt^2}h(t)\neq0$, there exists
a solution $\{q(t),p(t)\}$ of (\ref{pain7}), where $q(t)$ solves
 (\ref{pain1}).
An explicit correspondence between three sets $\{q(t),q'(t)\}$,
$\{q(t),p(t)\}$ and $\{h(t),h'(t),h''(t)\}$
is given by birational transformations,
which can be found in \cite{Ok87}.

The group of Backlund transformations of ${\bf P}_{VI}$ is
isomorphic to the affine Weyl group of the type $F_4$: $W_a(F_4)$.
It contains
the following transformations of parameters (only five of them
are independent)
\beq
w_1: b_1\leftrightarrow b_2,\quad w_2:
b_2\leftrightarrow b_3,\quad w_3: b_3\leftrightarrow b_4,
\quad w_4: b_3\to -b_3,\> b_4\to-b_4,\label{pain12}
\eeq

\beq
x^1: \kappa_0\leftrightarrow\kappa_1,\quad
x^2: \kappa_0\leftrightarrow\kappa_\infty,\quad
x^3: \kappa_0\leftrightarrow\theta\label{pain13}
\eeq
and the parallel transformation
\beq
l_3: {\bf b}\equiv(b_1,b_2,b_3,b_4)
\to{\bf b^+}\equiv(b_1,b_2,b_3+1,b_4).\label{pain14}
\eeq

The auxiliary function $h_+(t)$ corresponding to parameters
${\bf b^+}=l_3({\bf b})$
is given in \cite{Ok87}
\beq
h_+(t)=h(t)-q(q-1)p+(b_1+b_4)q
-\frac{1}{2}(b_1+b_2+b_4).\label{pain15}
\eeq
Following \cite{Ok87} one can calculate $h_+(t)$ in terms
of $h(t)$ and
its first and second derivatives
\beq
h_+(t)=\frac{t(t-1)h''(t)+2h(t)[b_3(b_3+1)+
h'(t)]+b_3(1-2t)h'(t)-b_1b_2b_4}
{2(h'(t)+b_3^2)}\label{pain16}
\eeq
and vice versa
\beq
h(t)=\frac{t(t-1)h_+''(t)+2h_+(t)[b_3(b_3+1)+h_+'(t)]-
(b_3+1)(1-2t)h_+'(t)+b_1b_2b_4}
{2(h_+'(t)+(b_3+1)^2)}.\label{pain17}
\eeq
For each solution of the ${\bf P}_{ VI}$ equation one can
introduce a corresponding
tau-function via
\beq
H(t,q(t),p(t);{\bf b})=\frac{d}{dt}
\log\mathrm{T}(t,{\bf b}).\label{pain18}
\eeq
Obviously tau-functions are defined up to an arbitrary
normalization factor.

Following \cite{Ok87} let us introduce a family of tau-functions
$\mathrm{T}_m(t)$
\beq
\mathrm{T}_m(t)=\exp{\left\{\int dt\, H(t,q(t),p(t);
{\bf b_m})\right\}},\label{pain19}
\eeq
where
\beq
{\bf b_m}\equiv l_3^m({\bf b})=(b_1,b_2,b_3+m,b_4),
\quad m\in {\bZ}.\label{pain20}
\eeq

As shown in \cite{Ok87} they satisfy the second order
Toda-type equation
\beq
\frac{d}{dt}\left[t(t-1)\frac{d}{dt}\log\mathrm{T}_m(t)\right]
+(b_1+b_3+m)(b_3+b_4+m)=c(m)
\frac{\mathrm{T}_{m+1}(t)\mathrm{T}_{m-1}(t)}
{\mathrm{T}^2_m(t)},\label{pain21}
\eeq
where $c(m)$ is a nonzero constant.

\section{Elliptic functions and useful identities}

In this section we list all definitions and properties
of elliptic functions used later in the text. Following
\cite{WW} we will use the standard theta functions
$\theta_i(x|\tau)$ with quasi-periods $\pi$ and $\pi\tau$.

The elliptic modulus  $k$ and its complement $k'$ are
defined in a standard way by
\beq
k=\frac{\theta_2(0|\tau)^2}{\theta_3(0|\tau)^2},\quad
k'=\frac{\theta_4(0|\tau)^2}{\theta_3(0|\tau)^2},\quad
k^2+k'^2=1.\label{pain22}
\eeq
It will be more convenient to use the parameter $t=k^2$
as the second argument of
elliptic functions and hereafter we will follow this notation
(except for the theta-functions), i.e.
\beq
\tau=i\frac{K'(t)}{K(t)},\quad q=e^{i\pi\tau},
\quad K(t)=\frac{\pi}{2}\theta_3^2(0|\tau),\quad K'(t)=K(1-t)
\label{pain23}
\eeq
where
${ K}(t)$ and ${ K}'(t)$ are the complete elliptic integrals
of the first kind of the
parameters $t$ and $1-t$.

We introduce Jacobi elliptic functions
\beq
\mbox{sn}(u,t)=\frac{1}{k^{1/2}}
\frac{\theta_1(v|\tau)}{\theta_4(v|\tau)},\quad
\mbox{cn}(u,t)=\frac{(k')^{1/2}}{k^{1/2}}
\frac{\theta_2(v|\tau)}{\theta_4(v|\tau)},\quad
\mbox{dn}(u,t)={(k')^{1/2}}
\frac{\theta_3(v|\tau)}{\theta_4(v|\tau)},\label{pain24}
\eeq
where
\beq
u=\frac{2K(t)}{\pi}v\label{pain24a}
\eeq
and define the fundamental elliptic integral of the
second kind \cite{WW} by
\beq
\mathcal {E}(u,t)=\int_0^u\mbox{dn}^2(x,t)dx.\label{pain25}
\eeq
It satisfies
\beq
\mathcal{E}(K(t),t)=E(t),\label{pain26}
\eeq
where $E(t)$ is the complete elliptic integral of the second kind.

Using (\ref{pain23}) and identities  for complete elliptic
integrals of the first and second kind
\beq
\p_tK(t)=\frac{E(t)}{2t(1-t)}-\frac{K(t)}{2t},\quad
\p_tE(t)=\frac{E(t)-K(t)}{2t}\label{pain27}
\eeq
\beq
E(t)K(1-t)+E(1-t)K(t)-K(1-t)K(t)=\pi/2,\label{pain28}
\eeq
one can obtain
\beq
\p_t\tau(t)=\frac{i\pi}{4t(t-1)K^2(t)}.\label{pain29}
\eeq

We also need the derivatives of Jacobi elliptic functions
with respect
to the parameter $t$.
Differentiating the formula
\beq
u=\sideset{}{}\int_0^{\mbox{\scriptsize sn}(u,t)}
\frac{dx}{\sqrt{(1-x^2)(1-tx^2)}} \label{pain30}
\eeq
with respect to $t$ and calculating the remaining integral we get
\beq
\frac{d}{d t}\mbox{sn}(u,t)=
-\frac{\mbox{sn}(u,t)\mbox{cn}^2(u,t)}{2(t-1)}+
\frac{\mbox{cn}(u,t)\mbox{dn}(u,t)}{2t(t-1)}\bigl[u(t-1)+
\mathcal{E}(u,t)\bigr].\label{pain31}
\eeq
From (\ref{pain31}) it is easy to obtain
\beq
\frac{d}{d t}\mbox{cn}(u,t)=
\frac{\mbox{sn}^2(u,t)\mbox{cn}(u,t)}{2(t-1)}-
\frac{\mbox{sn}(u,t)\mbox{dn}(u,t)}{2t(t-1)}\bigl[u(t-1)+
\mathcal{E}(u,t)\bigr]\label{pain32}
\eeq
and
\beq
\frac{d}{d t}\mbox{dn}(u,t)=
\frac{\mbox{sn}^2(u,t)\mbox{dn}(u,t)}{2(t-1)}-
\frac{\mbox{sn}(u,t)\mbox{cn}(u,t)}{2(t-1)}\bigl[u(t-1)+
\mathcal{E}(u,t)\bigr].\label{pain33}
\eeq

Integrating the well known formula between
a logarithmic derivative of $\theta_4(x|\tau)$ and $\mathcal E(u,t)$
we obtain
\beq
\theta_4(x|\tau)=\theta_4(0|\tau)
\exp\left\{-\frac{2x^2}{\pi^2}E(t)K(t)+
\sideset{}{}\int_0^{2xK(t)/\pi}\mathcal{E}(y,t)dy\right\}.
\label{pain34}
\eeq
We can use (\ref{pain34}) for calculation of the derivatives
$\theta_4(x|\tau)'_x$ and
$\theta_4(x|\tau)''_x$.

Finally, combining (\ref{pain34}) with the differential
equation satisfied by
theta-functions
\beq
\frac{4}{i\pi}\frac{\p}{\p\tau}\theta_i(u|\tau)+
\frac{\p^2}{\p u^2}\theta_i(u|\tau)=0,\quad i=1,2,3,4,
\label{pain35}
\eeq
one can calculate the derivative ${\p_\tau}\theta_4(x|\tau)$
 in terms of theta-functions and $\mathcal{E}(u,t)$.

We notice that formula (\ref{pain34}) is convenient
for expansion of $\theta_4(x|\tau)$ in a series in $x$
up to any required
order. Say,
\beq
\theta_4(x|\tau)=
\theta_4(0|\tau)\exp\left\{-\frac{2x^2}{\pi^2}E(t)K(t)\right\}
\left[1+2\frac{x^2 K^2(t)}{\pi^2}+
\frac{2(3-2t)}{3}\frac{x^4 K^4(t)}{\pi^4}+O(x^6)\right].
\label{pain36}
\eeq

\section{Picard solution and tau-functions}

The restriction on parameters (\ref{pain2}) for the
Picard solution (\ref{pain3})
can be rewritten in terms of parameters $b_i$ (\ref{pain6}) as
\beq
b_1=b_2=0,\quad b_3=b_4=-1/2.\label{pain37}
\eeq
It is convenient  to fix a particular branch of the Picard
solution (\ref{pain3}) by choosing two linearly independent
solutions
of  equation (\ref{pain4})
\beq
\omega_1(t)=\frac{\pi}{2}\phantom{|}_2F_1(\frac{1}{2},
\frac{1}{2};1;t)={ K}(t),\quad
\omega_2(t)=i\frac{\pi}{2}\phantom{|}_2F_1(\frac{1}{2},
\frac{1}{2};1;1-t)=i{ K}'(t)=i
K(1-t),
\label{pain38}
\eeq
where
${ K}(t)$ and ${ K}'(t)$ are the complete elliptic integrals
of the first kind as defined in the previous section.

With a choice of half-periods (\ref{pain38}) for the
Weierstrass function
the expressions for the invariants  $e_1, e_2, e_3$
take the following form
\beq
e_1=1-\frac{t+1}{3},
\quad e_2=t-\frac{t+1}{3},\quad e_3=-\frac{t+1}{3}.
\label{pain39}
\eeq
After some simple calculations
we can rewrite the Picard solution of ${\bf P}_{VI}$ as
\beq
q_0(t)=\frac{1}{\mbox{sn}^2(c_1{K}(t)+ic_2{ K}'(t),t)}
,\label{pain40}
\eeq
where  $c_1$ and $c_2$ are the same parameters as in
(\ref{pain3}).

Let us make a change of variables
\beq
c_1=\frac{2x}{\pi}+1,\quad c_2=\frac{2y}{\pi}+1,\label{pain41}
\eeq
where $x$, $y$ are new parameters. A reason for this is that
the resulting tau-functions look simpler in $x$ and $y$.

With a substitution (\ref{pain41}) formula (\ref{pain40})
takes the form
\beq
q_0(t)=t\,\frac{\mbox{cn}^2(z,t)}{\mbox{dn}^2(z,t)},
\label{pain42}
\eeq
where we defined a new variable $z$
\beq
z=\frac{2{ K(t)}}{\pi}(x+\tau y).\label{pain43}
\eeq
The hamiltonian $H_0(t)$ for the choice of parameters
(\ref{pain37})
can be calculated from (\ref{pain7}-\ref{pain8}) in terms
of $q_0(t)$ and its first derivative
\beq
H_0(t)=\frac{t^2+(1-2t)q_0(t)}{4t(t-1)(q_0(t)-t)}+
\frac{t(t-1)q_0'(t)^2}{4q_0(t)(q_0(t)-1)(q_0(t)-t)}.
\label{pain44}
\eeq

Using formulas from section 3 one can explicitly calculate
the derivative $q_0'(t)$

\beq
\frac{d}{dt}q_0(t)=q_0'(t)=
\frac{1}{\mbox{dn}^2(z,t)}+
\frac{\mbox{sn}(z,t)\mbox{cn}(z,t)}{\mbox{dn}^3(z,t)}
\left[\frac{\pi}{2K(t)}
\log\theta_2(x+\tau y|\tau)]'_x+\frac{i y}{K(t)}
\right].\label{pain45}
\eeq
Substituting (\ref{pain45}) into (\ref{pain44}) we produce
 the following
expression for the function $H_0(t)$:
\beq
H_0(t)=-\frac{1}{4(t-1)}-
\frac{\mbox{cn}^2(z,t)}{4t(t-1)\,\mbox{sn}^2(z,t)}+
\frac{\mathscr{E}(x,y,t)^2}{4t(t-1)}, \label{pain46}
\eeq
where
\beq
\mathscr{E}(x,y,t)=
\frac{\pi}{2K(t)}[\log\theta_1(x+\tau y|\tau)]'_x+
\frac{iy}{K(t)}\label{pain47}
\eeq
and $z$ is defined by (\ref{pain43}).

To calculate the tau-function $\mathrm{T}_0(t)$ for the
Picard solution
(\ref{pain42}) we have to calculate the indefinite integral
of (\ref{pain46})
with respect to the variable $t$ which looks like
a hopeless problem.

Now we formulate the central result of this paper.
\begin{theorem}
The tau-function for the Picard solution (\ref{pain42})
is given by
\beq
\mathrm{T}_0(t)=\exp\{\int H_0(t)dt\}=c_0(x,y)
\,{q^{{y^2/{\pi^2}}}}{t^{-1/4}}\,
\frac{\theta_1(x+\tau y|\tau)}
{\theta_4(0|\tau)}, \label{pain48}
\eeq
where $c_0(x,y)$ is an integration constant.
\end{theorem}

{\bf Proof:}

First we rewrite $\mathrm{T}_0(t)$ as
\beq
\mathrm{T}_0(t)=
c_0(x,y)\, q^{{y^2}/{\pi^2}}\,\mbox{sn}(z,t)
\frac{\theta_4(x+\tau y|\tau)}
{\theta_4(0|\tau)}.\label{pain49}
\eeq
The proof is straightforward and reduces to
differentiations.
We shall do this in a few steps.
Taking logarithmic derivative of (\ref{pain49})
we obtain
\bea
&&\frac{\p}{\p t}\log \mathrm{T}_0(t)=
i\frac{y^2}{\pi}\p_t\tau(t)+
\frac{\mbox{cn}(z,t)\mbox{dn}(z,t)}{\mbox{sn}(z,t)}
\left[\frac{2(x+\tau y)}{\pi}\p_tK(t)+
\frac{2yK(t)}{\pi}\p_t\tau(t)\right]
+
\nonumber\\
&&+\frac{\p_t\,\mbox{sn}(u,t)}{\mbox{sn}(u,t)}\biggr|_{u=z}+
\left[\left(y\left[\log\theta_4(u|\tau)\right]'_u+
\frac{\p_\tau\theta_4(u|\tau)}
{\theta_4(u|\tau)}\right)\biggr|_{u=x+\tau y}-
\frac{{\p_\tau}\theta_4(0|\tau)}{\theta_4(0|\tau)}\right]
\p_t\tau(t).\label{pain50}
\eea

Using (\ref{pain27}-\ref{pain29}) one can evaluate the
first two terms in (\ref{pain50}).
The derivative ${\p_t}\mbox{sn}(u,t)$ was calculated in
(\ref{pain31}).
Differentiating  (\ref{pain34}) twice and using the
equation (\ref{pain35})
one can evaluate
all derivatives in the last term of (\ref{pain50}).
Combining all contributions and using
a simple formula
\beq
\log\theta_1(x+\tau y|\tau)]'_x=
\log\theta_4(x+\tau y|\tau)]'_x+
\frac{2K(t)}{\pi}
\frac{\mbox{cn}(z,t)\mbox{dn}(z,t)}{\mbox{sn}(z,t)},
\label{pain51}
\eeq
we obtain after simplifications
  the expression (\ref{pain46}) for $H_0(t)$.
\hfill\qed
\vspace{0.5cm}

It is quite remarkable that the indefinite integral of
the function (\ref{pain46})  gives such a simple answer
(\ref{pain48}).
We were able to produce this expression by expanding
(\ref{pain46})
in a series in $x$ at $y=0$ and integrating term by term.
Analyzing  the resulting series we compared  it to
the expansion
of $\theta_1(x|\tau)$ in $x$ which is similar to
the expansion (\ref{pain36}).
This allowed us to arrive  at the final answer (\ref{pain48}).

To solve the equation (\ref{pain21}) for $\mathrm{T}_m(t)$
we need to calculate
the second tau-function $\mathrm{T}_1(t)$
corresponding to the solution with
parameters
\beq
{\bf b_1}=(0,0,1/2,-1/2).\label{pain52}
\eeq
First we have from (\ref{pain10})
\beq
h_0(t)=t(t-1)H_0(t)+t/4-1/8\label{pain53}
\eeq
for the Picard solution with ${\bf b_0}=(0,0,-1/2,-1/2)$ and
\beq
h_1(t)=t(t-1)H_1(t)-t/4+1/8\label{pain54}
\eeq
for the solution with parameters (\ref{pain52}).

Now $h_1(t)$ is obtained using the birational canonical
transformation
 (\ref{pain16}) with parameters (\ref{pain37}) and
$h(t)$ replaced with $h_0(t)$. Using
(\ref{pain53}-\ref{pain54})
one can arrive at the following answer:
\beq
H_1(t)=-\frac{\mbox{sn}^2(z,t)}{4\,\mbox{dn}^2(z,t)}+
\frac{1}{4t(t-1)}
\left[\frac{\pi}{2K(t)}[\log\theta_4(x+\tau y|\tau)]'_x+
\frac{iy}{K(t)}-
t\frac{\mbox{sn}(z,t)\mbox{cn}(z,t)}{\mbox{dn}(z,t)}
\right]^2.\label{pain55}
\eeq
In fact, the equation satisfied by $h_0(t)$ and $h_1(t)$
is the same
and it is easy to check that
\beq
h_1(t,x,y)=h_0(t,x+\frac{\pi}{2},y+\frac{\pi}{2})
\label{pain55a}
\eeq
and
\beq
H_1(t,x,y)=H_0(t,x+\frac{\pi}{2},y+\frac{\pi}{2})+
\frac{1}{4t}+\frac{1}{4(t-1)},\label{pain55b}
\eeq
where we show explicitly a dependence on the fixed
parameters $x$ and $y$.

Taking (\ref{pain55b}) into account it is easy to
integrate $H_1(t)$ using  formula (\ref{pain48}).

The answer is given by the following
\begin{theorem}
The $\tau$-function $\mathrm{T}_1(t)$ is given by
\beq
\mathrm{T}_1(t)=\exp\{\int H_1(t)dt\}=
c_1(x,y)\,{q^{{y^2/{\pi^2}}}}(1-t)^{1/4}\,
\frac{\theta_3(x+\tau y|\tau)}
{\theta_4(0|\tau)}, \label{pain56}
\eeq
where $c_1(x,y)$ is an arbitrary integration constant.
\end{theorem}
\vspace{0.3cm}

If we define  the sequence of tau-functions $\mathrm{T}_m(t)$
corresponding to Picard type solutions with parameters
\beq
{\bf b}_m=(0,0,-1/2+m,-1/2),\quad m\in \bZ\label{pain57}
\eeq
and two initial conditions
\beq
\mathrm{T}_0(t)={q^{{y^2/{\pi^2}}}}{t^{-1/4}}\,
\frac{\theta_1(x+\tau y|\tau)}
{\theta_4(0|\tau)},\quad
\mathrm{T}_1(t)={q^{{y^2/{\pi^2}}}}(1-t)^{1/4}\,
\frac{\theta_3(x+\tau y|\tau)}
{\theta_4(0|\tau)},
\eeq
then other tau-functions $\mathrm{T}_m(t)$, for
$m>1$ or $m<0$ can be
calculated from the difference-differential equation
\beq
\frac{d}{dt}\left[t(t-1)
\frac{d}{dt}\log\mathrm{T}_m(t)\right]
+\left(m-\frac{1}{2}\right)^2=c(m)
\frac{\mathrm{T}_{m+1}(t)\mathrm{T}_{m-1}(t)}
{\mathrm{T}^2_m(t)},\label{pain58}
\eeq
where $c(m)$ is determined by a normalization of tau-functions.

Note that the expressions for $\mathrm{T}_m(t)$, $m\neq0,1$
will be more complicated
and involve explicitly the function $\mathcal{E}(u,t)$
defined in (\ref{pain25}).
We will not calculate them here.

\section{Conclusion}

In this paper we constructed a sequence of tau-functions for the
Picard type solutions of Painlev\'e VI equation with parameters
(\ref{pain57}). In fact, starting with the tau-function for the
Picard solution
(\ref{pain48}) one can write many different
birational canonical transformations and calculate corresponding
sequences of tau-functions.
We have successfully  applied this approach to
sum up some infinite form factor expansions for the 2D Ising model.
All details will be given in  forthcoming publications.

\section*{Acknowledgments}
I would like to thank M.T. Batchelor,  J. De Gier,
S.M. Sergeev for useful remarks
and V.V. Bazhanov,  A.J. Guttmann for stimulating discussions and
careful reading of
the manuscript.
This work has been supported by the Australian Research Council.

\end{document}